\documentclass[12pt]{amsart}
\let\displaystyle\textstyle
\usepackage{amsmath,amstext,amsfonts,amssymb,mathrsfs,mathtools,setspace}
\usepackage[alphabetic]{amsrefs}
\usepackage[unicode]{hyperref}
\usepackage[utf8]{inputenc}
\usepackage[T1]{fontenc}
\usepackage{lmodern}
\hypersetup{colorlinks=false,%
linkbordercolor={0.93 0.71 0.13},%
citebordercolor={1 0 .5},%
urlbordercolor={0.27 0.39 0.68}}

\newcommand{\KAKEYA}{\SCR K}
\newcommand{\suchthat}{\mid}

\newcommand{\RR}{\bb R}

\newcommand{\PROB}{\mathbb{P}}

\newcommand{\m}[1]{\(\smash{#1}\)}
\newcommand{\algn}[1]{\begin{align*} #1\end{align*}}

\newcommand{\bb}[1]{\mathbb{#1}}
\newcommand{\QQ}{\bb{Q}}
\newcommand{\FF}{\bb{F}}
\newcommand{\ZZ}{\bb{Z}}

\newcommand{\SCR}{\mathscr}

\newcommand{\PADICVALUATION}{\smash{v_{p}}}

\newcommand{\TADICVALUATION}{\smash{v_{t}}}

\newcommand{\GEQSLANT}{\geqslant}
\newcommand{\LEQSLANT}{\leqslant}

\newcommand{\maxideal}{\mathfrak m}

\newcommand{\Hypergeometric}[5]%
{\smash{\null_{#1}\mathrm{F}_{#2}\big(%
\smash{\genfrac{}{}{0pt}{2}{#3}{#4}};{#5}\big)}}

\newtheorem{theorem}{Theorem}
\newtheorem{definition}[theorem]{Definition}
\newtheorem{lemma}[theorem]{Lemma}

\newtheorem{conjecture}{Conjecture}

\theoremstyle{remark}
\newtheorem{remark}[theorem]{Remark}

\makeatletter
\newcommand{\dashover}[2][\mathop]{#1{\mathpalette\df@over{{\dashfill}{#2}}}}
\newcommand{\fillover}[2][\mathop]{#1{\mathpalette\df@over{{\solidfill}{#2}}}}
\newcommand{\df@over}[2]{\df@@over#1#2}
\newcommand\df@@over[3]{%
  \vbox{
    \offinterlineskip
    \ialign{##\cr
      #2{#1}\cr
      \noalign{\kern1pt}
      $\m@th#1#3$\cr
    }
  }%
}
\newcommand{\dashfill}[1]{%
  \kern-.5pt
  \xleaders\hbox{\kern.5pt\vrule height.4pt width \dash@width{#1}\kern.5pt}\hfill
  \kern-.5pt
}
\newcommand{\dash@width}[1]{%
  \ifx#1\displaystyle
    2pt
  \else
    \ifx#1\textstyle
      1.5pt
    \else
      \ifx#1\scriptstyle
        1.25pt
      \else
        \ifx#1\scriptscriptstyle
          1pt
        \fi
      \fi
    \fi
  \fi
}
\newcommand{\solidfill}[1]{\leaders\hrule\hfill}
\makeatother

    {\tiny (\begin{smallmatrix}}%
    {\end{smallmatrix})}
    {\footnotesize (\begin{smallmatrix}}%
    {\end{smallmatrix})}
\newcommand{\BLACKSQUARE}{\hspace{\stretch{1}}\m{\mathbin{\vcenter{%
    \hbox{\rule{.9ex}{.9ex}}}}}}
\renewenvironment{proof}[1]%
    {\emph{{#1}}}%
    {\BLACKSQUARE}

\makeatletter
\def\@settitle{\begin{center}%
  \baselineskip14\p@\relax\normalfont\LARGE
    \uppercasenonmath\@title
  \@title\end{center}%
}
\def\do#1{\@ifdefinable#1{\newdimen#1}}
\do\@TempWidthOne
\do\@TempWidthTwo
\newcommand*\WidthInMathUnitsOf[2]{%
    \settowidth\@TempWidthOne{$\m@th #1#2$}%
    \settowidth\@TempWidthTwo{$\m@th #1\mspace{1mu}$}%
    \strip@pt \dimexpr \@TempWidthOne*\p@/\@TempWidthTwo \relax
}
\newcommand*\SetToWidthInMathUnits[3]{%
    \settowidth\@TempWidthOne{$\m@th #2#3$}%
    \settowidth\@TempWidthTwo{$\m@th #2\mspace{1mu}$}%
    \edef #1{%
        \strip@pt \dimexpr \@TempWidthOne*\p@/\@TempWidthTwo \relax
    }%
}
\newcommand{\llrrparen}[1]{
  (\!({#1})\!)}

\newcommand*\ExtractHexDigitOfFamilyOfSymbol[1]{%
    \expandafter\@ExtractHexDigitOfFamily \meaning #1%
}
\@ifdefinable\@ExtractHexDigitOfFamily{
    \def\@ExtractHexDigitOfFamily#1"#2#3#4#5{#3}
}
\makeatother

\newcommand*{\ReportValueOfOneEm}[1]{%
    &\text{in #1 size,}%
        &\quad\texttt{1em}%
            &=\texttt{\the \fontdimen 6 \csname #1font\endcsname 2}%
        &&\qquad\text{(from \texttt{\fontname \csname #1font\endcsname 2})}%
}
\newcommand*{\ReportSizeOfFontOfSymbol}[2]{%
    &\texttt{\pdffontsize\csname #2font\endcsname
                "\ExtractHexDigitOfFamilyOfSymbol{#1}}%
        &&\quad\text{in #2 size}%
        &\qquad(\texttt{\string\fontname}%
            &=\texttt{\fontname\csname #2font\endcsname
                "\ExtractHexDigitOfFamilyOfSymbol{#1}})%
}

\title%
[The \m{p}-adic Kakeya conjecture]%
{\textrm{The \m{p}-adic Kakeya conjecture}}
\author{Bodan Arsovski}
\date{}
\address{School of Mathematics and Statistics, University of Sheffield}
\email{\href{mailto:bodan.arsovski@outlook.com}{bodan.arsovski@outlook.com}}
\urladdr{\href{https://bodanarsovski.github.io/}{bodanarsovski.github.io}}


\begin{document}
\maketitle

\begin{abstract}
We prove that all bounded subsets of \m{\QQ_p^n} containing a line segment of unit length in every direction have Hausdorff and Minkowski dimension \m{n}. This is the analogue of the classical Kakeya conjecture with \m{\RR} replaced by \m{\QQ_p}.
\end{abstract}

\section{Introduction}
\label{secInt}

In 1917, Kakeya posed the Kakeya needle problem, asking about the minimum area of a region in the plane in which a needle of unit length can be rotated around by \m{360^\circ}. The problem reduces to asking about the minimum area of a \emph{Kakeya set} in \m{\RR^2}---that is, a bounded subset of \m{\RR^2} containing a line segment of unit length in every direction.  Besicovitch~\cite{Bes28} showed that, in a certain sense, the answer is ``arbitrarily small'', by constructing a Kakeya set in \m{\RR^2} of Lebesgue measure zero. On the other hand, Davies~\cite{Dav71} showed that, in a different sense, the answer is ``large'': a Kakeya set in \m{\RR^2} must have Minkowski dimension \m{2}. The construction of~\cite{Bes28} extends immediately to \m{\RR^n} with \m{n\GEQSLANT 3}.
 The higher-dimensional analogue of the result of~\cite{Dav71} is much more difficult: it is the Minkowski dimension version of the notorious Kakeya conjecture, which is one of the most important open problems in geometric measure theory, and analysis in general.

\begin{conjecture}[Kakeya]\label{Kak1}
Let \m{n} be a positive integer. All Kakeya sets in \m{\RR^n} have Hausdorff and Minkowski dimension \m{n}.
\end{conjecture}

The Kakeya conjecture has deep connections with harmonic analysis, PDEs, and combinatorics among other fields, and it is open for \m{n\GEQSLANT 3}.
For \m{n\in\{3,4\}}, the state of the art are the results of Katz--Zahl~\cites{KZ19,KZ21} that all Kakeya sets in \m{\RR^3}  have Hausdorff and Minkowski dimension at least \m{5/2+\varepsilon} (for some \m{\varepsilon>0}), and all Kakeya sets in \m{\RR^4} have Hausdorff and Minkowski dimension at least \m{3.059} (see also Wolff~\cite{Wol95}, Katz--Łaba--Tao~\cite{KLT00}, and Łaba--Tao~\cite{LT01} for previous results). For \m{n\GEQSLANT 5}, the state of the art is the result of Katz--Tao~\cite{KT02} that all Kakeya sets in \m{\RR^n} have Hausdorff dimension at least \m{(2-\sqrt{2})(n-4)+3} and Minkowski dimension at least \m{(n-1)/\alpha+1}, where \m{\alpha \approx 1.675} is the largest root of the polynomial \m{z^3-4z+2}.

As a possible approach to the Kakeya conjecture, Wolff~\cite{W99} suggested an analogous conjecture over finite fields, which was then proved by Dvir~\cite{D09}. As noted by Ellenberg--Oberlin--Tao~\cite{EOT10}, the analogy between the classical and  finite field Kakeya conjectures 
 is imperfect for two reasons. Firstly, there is no natural non-discrete notion of distance in finite vector spaces, and the discrete distance does not have multiple scales. By contrast, the Euclidean distance in \m{\RR^n} has multiple scales.
 Secondly, the result of~\cite{D09} is \emph{too strong:} all Kakeya sets in \m{\FF_q^n} (where  \m{\FF_q} is a finite field) have, in a certain sense, positive measure. In particular, there is no variant in \m{\FF_q^n} of the construction of Besicovitch of a Kakeya set of measure zero. 
 
In light of these two shortcomings, Ellenberg--Oberlin--Tao~\cite{EOT10} suggested analogues of the Kakeya conjecture over topological fields whose metrics capture the multiple scales property of the Euclidean metric of \m{\RR^n}. The only (non-discrete) locally compact Hausdorff topological fields are the finite extensions of \m{\RR}, \m{\QQ_p}, and \m{\FF_q\llrrparen{t}} (for varying primes \m{p} and prime powers \m{q}). Moreover, the Kakeya conjecture over \m{K \in \{\RR,\QQ_p,\FF_q\llrrparen{t}\}} implies the Kakeya conjecture over any finite extension \m{L/K}, simply by treating \m{L^n} as a vector space over \m{K}. Therefore, other than the classical setting over \m{\RR}, the two other natural settings for the Kakeya conjecture (with interesting metrics that have multiple scales) are over \m{\QQ_p} and \m{\FF_q\llrrparen{t}}.
  Both shortcomings are resolved in these settings. Firstly, \m{\QQ_p^n} and \m{\FF_q\llrrparen{t}^n} have metrics that have multiple scales, just like the metric of \m{\RR^n}. Secondly, there are variants of the Besicovitch construction in both settings: Dummit--Hablicsek~\cite{DH11} constructed a Kakeya set in \m{\FF_q\llrrparen{t}^n} of measure zero, Hickman--Wright~\cite{HW18} constructed a Kakeya set in %
   \m{\QQ_p^n} of measure zero, and Fraser~\cite{Fra16} constructed more generally a Kakeya set in \m{L^n} of measure zero, for any finite extension \m{L/K} of \m{K \in \{\QQ_p,\FF_q\llrrparen{t}\}}. Moreover, Dummit--Hablicsek~\cite{DH11}  proved the %
   analogues of the result of Davies~\cite{Dav71} in \m{\QQ_p^2} and \m{\FF_q\llrrparen{t}^2}. %
    All of these theorems suggest that the Kakeya conjectures over \m{\QQ_p} and \m{\FF_q\llrrparen{t}} could be good models for the classical Kakeya conjecture over \m{\RR}.
 
 In this article we prove the Kakeya conjecture over \m{\QQ_p}: %

\begin{theorem}\label{mth}
Let \m{p} be a prime number and \m{n} a positive integer. All Kakeya sets in \m{\QQ_p^n} have Hausdorff and Minkowski dimension \m{n}.
\end{theorem}

\begin{remark}
Let \m{p} and \m{n} be as in theorem~\ref{mth}, and let \m{K/\QQ_p} be a finite extension.
Theorem~\ref{mth} then implies that all Kakeya sets in \m{K^n} have Hausdorff and Minkowski dimension \m{n}, i.e.~the Kakeya conjecture over all \m{p}-adic local fields (simply by treating \m{K^n} as a vector space over \m{\QQ_p}). %
\end{remark}

\begin{remark}
 The field of \m{p}-adic numbers \m{\QQ_p} is the completion of \m{\QQ} with respect to the \m{p}-adic norm \m{|x|_p = p^{-\PADICVALUATION(x)}}, and, in fact, the only completions of \m{\QQ} are \m{\RR} and \m{\QQ_p} (for varying \m{p}). %
  To give some intuition behind its topology, %
  the closed unit ball \m{\ZZ_p} of \m{\QQ_p} is homeomorphic to the Cantor set \m{\SCR C\subseteq [0,1]}. There is a well developed analogue of classical analysis over the \m{p}-adic numbers called {\m{p}-adic analysis}, which has historically been of particular interest in number theory.
\end{remark}

The proof of theorem~\ref{mth} is inspired by Dvir's famous proof of the finite field Kakeya conjecture~\cite{D09}, and the method can be summarized as ``the polynomial method over discrete valuation rings''. This discrete valuation polynomial method seems flexible and extensible, and we hope that it can be used to adapt other theorems with known polynomial method proofs, from the finite field setting over \m{\FF_q}, to the \m{p}-adic setting over \m{K/\QQ_p} (and in particular over the finite rings \m{\ZZ/p^k\ZZ}).
We prioritize brevity and clarity of exposition, so the implied constants are not sharp and can be improved by more careful analysis\footnote{In particular, in an earlier version of this manuscript we gave a different proof inspired by ideas of Dhar--Dvir~\cite{DD}, that results in slightly better constants. To keep the exposition brief and simple, we omit that proof in favor of the discrete valuation polynomial method proof, which feels more conceptual.}.

\section{Proof}
\label{set}

Let \m{p} be a prime number, \m{n} and \m{k} positive integers, and \m{q=p^k}. 
 Let \m{\FF=\FF_p}, and \m{R=\ZZ/q\ZZ}. 
  Let \m{\QQ_p} denote the \m{p}-adic numbers, %
 and \m{\ZZ_p} the \m{p}-adic integers. That is, \m{\QQ_p} is the completion of \m{\QQ} with respect to the norm \m{|x|_p = p^{-\PADICVALUATION(x)}}, and \m{\ZZ_p} is the closed unit ball in \m{\QQ_p} centered at the origin %
 \m{\ZZ_p = \{z \in \QQ_p \suchthat |z|_p \LEQSLANT 1\}\subseteq \QQ_p}.  
  We equip  \m{\QQ_p^n} with the maximum product norm \algn{|(x_1,\ldots,x_n)|_p = \max_{i \in \{1,\ldots,n\}}|x_i|_p, \text{ for } (x_1,\ldots,x_n) \in \QQ_p^n,} and by the Hausdorff and Minkowski dimensions of a subset of \m{\QQ_p^n} we of course mean the Hausdorff and Minkowski dimensions with respect to this norm. Let \m{\zeta} be a primitive \m{q^\text{th}} root of unity in the algebraic closure of \m{\QQ_p}. Let
\m{T=\ZZ_p(\zeta)[t]}, and let \m{\overline T=\FF[t]} be the reduction of \m{T} modulo the maximal ideal \m{\maxideal = (\zeta-1)} of \m{\ZZ_p[\zeta]}. Thus \m{\overline T} is a subset of the discrete valuation ring \m{\overline T_{(t)}}, whose field of fractions is \m{\FF(t)}, and whose valuation we denote by \m{\TADICVALUATION}. Let 
\m{C = \{1,\ldots,(1+t)^{q-1}\} \subseteq T}, and let \m{\overline C \subseteq \overline T} be the reduction of \m{C} modulo \m{\maxideal}.
\begin{definition}\label{kakdef}
Let \m{\varepsilon,\delta\in[0,1]}. Let us normalize the measures of the closed unit balls \m{\ZZ_p\subseteq \QQ_p} and \m{\ZZ_p^n\subseteq \QQ_p^n} to be equal to \m{1}.

An \m{(\varepsilon,\delta)}-Kakeya set in \m{\QQ_p^n} is a bounded subset of \m{\QQ_p^n} containing at least \m{\varepsilon} of a line segment of unit length in at least \m{\delta} of directions---that is, a bounded subset \m{\KAKEYA\subseteq \QQ_p^n} with the following property: there is a subset \m{D\subseteq \ZZ_p^n} of measure at least \m{\delta} such that, for all \m{x\in D}, there are \m{b_x\in \QQ_p^n} and a subset \m{\Lambda_x \subseteq \ZZ_p} of measure at least \m{\varepsilon} such that \m{b_x+\lambda x\in \KAKEYA} for all \m{\lambda \in \Lambda_x}.\footnote{The definition of a Kakeya set in \m{\QQ_p^n} is slightly different in~\cites{EOT10,HW18,DD} (only directions in \m{\bb P^{n-1}(\ZZ_p)} are considered), but it is equivalent to our definition.}
\end{definition}

We obtain theorem~\ref{mth} as a corollary to the following theorem, which gives a lower bound on the number of closed balls of radii \m{p^{-k}} needed to cover an \m{(\varepsilon,\delta)}-Kakeya set in the closed unit ball \m{\ZZ_p^n} of \m{\QQ_p^n}.

\begin{theorem}\label{qth}
Let \m{p} be a prime number, and \m{n} and \m{k} positive integers. Let \m{\varepsilon,\delta \in [0,1]}. 
An \m{(\varepsilon,\delta)}-Kakeya set %
  in \m{\ZZ_p^n}  
 cannot be covered by fewer than \algn{\binom{\varepsilon \delta q/pnk+n-1}{n}} closed balls of radii \m{p^{-k}}.
\end{theorem}

\begin{proof}{Proof that theorem~\ref{qth} implies theorem~\ref{mth}.}
Suppose that \algn{\varepsilon,\delta \gtrsim_{p,n} \frac{1}{\log_p(1/p^{-k})} = \frac1k.} An \m{(\varepsilon,\delta)}-Kakeya set \m{\KAKEYA\subseteq \QQ_p^n} is bounded, and therefore contained in finitely many disjoint closed unit balls. Suppose that \m{\KAKEYA} is covered by \m{N(\KAKEYA,p^{-k})} closed balls of radii \m{p^{-k}}. By translating and superimposing these closed unit balls onto the closed unit ball \m{\ZZ_p^n}, we end up with an \m{(\varepsilon,\delta)}-Kakeya set \m{\tilde \KAKEYA\subseteq \ZZ_p^n} that is covered by \m{N(\KAKEYA,p^{-k})} closed balls of radii \m{p^{-k}}, implying, by theorem~\ref{qth}, that \algn{N(\KAKEYA,p^{-k})\GEQSLANT  \frac{(\varepsilon \delta)^n p^{kn}}{(pnk)^n n!} \gtrsim_{p,n} \frac{p^{kn}}{k^{3n}}.} The standard dyadic method of Bourgain~\cite{Bou99} can be  used to conclude that the Hausdorff dimension of \m{\KAKEYA} is \m{n} (and therefore so is the Minkowski dimension).
\end{proof}

\begin{lemma}[Discrete valuation Schwartz--Zippel]\label{vanlem}
Let the coefficient of \m{z_1^{m_1}\cdots z_n^{m_n}} in \m{f \in \overline T[z_1,\ldots,z_n]} be \m{c\ne 0}, where \m{m_i \in \{0,\ldots,q\}} for all \m{i\in\{1,\ldots,n\}}, let the coefficient of any monomial that is larger than \m{z_1^{m_1}\cdots z_n^{m_n}} in the lexicographic order be zero, and let \m{\nu\in(0,1]}.
The number of \m{s \in \overline C^n} such that \m{\TADICVALUATION\left(f(s)\right) \GEQSLANT \TADICVALUATION\left(c\right)+\nu n q} is at most \m{pkq^{n-1} (m_1+\cdots+m_n)/\nu}.
\end{lemma}

\begin{proof}{Proof.}
We use induction. Let \m{n=1}, and suppose, on the contrary, that there is a subset \m{S\subseteq \overline C} such that \m{|S|>pkm_1/\nu}, on which \m{f} has \m{\TADICVALUATION}-valuation at least \m{\TADICVALUATION\left(c\right)+\nu q}. Let \m{S = \{(1+t)^{l} \suchthat l \in L\}}, so that \m{L\subseteq \{0,\ldots,q-1\}} and \m{|L|=|S|}. Then \m{L} has a subset \m{L_0} that contains no pair of distinct elements whose difference is divisible by \m{p^{-\lceil\log_p(k/\nu)\rceil}q}, and such that \m{|L_0| =m_1+1}. Let \m{S_0} be the subset of \m{S} that corresponds to \m{L_0}. By Lagrange interpolation (in the field \m{\FF(t)\supset \overline T}),
\algn{1 = \sum_{s \in S_0} \left(c\prod_{u \in S_0\backslash \{s\}} (s-u)\right)^{-1} f(s),}
so there exists some \m{s \in S_0} such that \algn{\TADICVALUATION\left(\prod_{u \in S_0\backslash \{s\}} (s-u)\right) \GEQSLANT \nu q.} %
 Let \m{L_0 = \{l_0,\ldots,l_{m_1}\}} with \m{l_0=s}. A theorem of Kummer %
  implies that \m{\binom{l}{w}} is a unit in \m{\FF} if and only if every \m{p}-adic digit of \m{w} is at most as large as the corresponding \m{p}-adic digit of \m{l}. So \algn{\TADICVALUATION\left((1+t)^{l_i}-(1+t)^{l_0}\right)=\TADICVALUATION\left((1+t)^{l_i-l_0}-1\right)=p^{\PADICVALUATION(l_i-l_0)},} and
\algn{\TADICVALUATION\left(\prod_{u \in S_0\backslash \{s\}} (s-u)\right) & = \sum_{i=1}^{m_1} \TADICVALUATION\left((1+t)^{l_i}-(1+t)^{l_0}\right) \\&=\sum_{i=1}^{m_1} p^{\PADICVALUATION(l_i-l_0)} \\ & \LEQSLANT \sum_{j=1}^{\smash{\lceil \log_pm_1\rceil}} p^{-\lceil\log_p(k/\nu)\rceil-j} q (p^j-p^{j-1}) \\ & < \nu q \lceil \log_pm_1\rceil/k \LEQSLANT \nu q,} a contradiction.
 That concludes the case \m{n=1}. For \m{n>1}, we write
\algn{f(z_1,\ldots,z_n) = z_1^{m_1} g(z_2,\ldots,z_n) + h(z_1,\ldots,z_n)} with the degree of \m{z_1} in \m{h} less than \m{m_1}. We evaluate \m{z_1,\ldots,z_n} on \m{\overline C} independently and uniformly at random. Let \m{F,G_{<}}, and \m{G_{\GEQSLANT}} be the events
\algn{F & \quad\quad \TADICVALUATION\left(f\right)\GEQSLANT \TADICVALUATION\left(c\right)+\nu n q,
\\ G_< & \quad\quad \TADICVALUATION\left(g\right)< \TADICVALUATION\left(c\right)+\nu(n-1) q,
\\ G_{\GEQSLANT} & \quad\quad \TADICVALUATION\left(g\right)\GEQSLANT \TADICVALUATION\left(c\right)+\nu(n-1) q.}
Then \algn{\PROB\left(F\right) & = \PROB\left(F\suchthat G_<\right) \PROB\left(G_<\right) + \PROB\left(F\suchthat G_{\GEQSLANT}\right)\PROB\left(G_{\GEQSLANT}\right)
\\ & \LEQSLANT \PROB\left(F\suchthat G_<\right) + \PROB\left(G_{\GEQSLANT}\right) \\& \LEQSLANT pkm_1/\nu q + pk(m_2+\cdots+m_n)/\nu q,} where the last line follows from the induction hypothesis, concluding the induction step, and with that the proof.
\end{proof}

\ \null\\[-19pt]

\begin{proof}{Proof of theorem~\ref{qth}.}
The image under the projection \m{\ZZ_p^n \to R^n} of an \m{(\varepsilon,\delta)}-Kakeya set in \m{\ZZ_p^n} is a subset of \m{R^n} %
 with the following property: there is a subset \m{D\subseteq R^n} of size at least \m{\delta |R|^n} such that, for all \m{x\in D}, there are \m{b_x\in R^n} and a subset \m{\Lambda_x \subseteq R} of size at least \m{\varepsilon q} such that \m{b_x+\lambda x\in S} for all \m{\lambda \in \Lambda_x}. %
  Our goal is to show that such a subset of \m{R^n} must have size at least \m{\binom{\varepsilon \delta q/pnk+n-1}{n}}.
Suppose, on the contrary, that there is such a subset \m{S \subseteq R^n} of size less than \m{\binom{\varepsilon \delta q/pnk+n-1}{n}}. In particular, \m{\varepsilon\delta>0}. Let \algn{S_C = \left\{(\zeta^{s_1},\ldots,\zeta^{s_n}) \in \ZZ_p[\zeta]^n \suchthat (s_1,\ldots,s_n)\in S\right\}.}
Then there is a  non-zero polynomial \m{f \in \QQ_p(\zeta)[z_1,\ldots,z_n]} of total degree less than \m{\varepsilon \delta q/pnk} that vanishes on \m{S_C}. 
After  multiplying \m{f} by \m{(1-\zeta)^m} for a suitable \m{m\in\ZZ}, we may assume that all coefficients of \m{f} have non-negative \m{\PADICVALUATION}-valuation and some coefficients of \m{f} have \m{\PADICVALUATION}-valuation equal to zero, implying that \m{f \in \ZZ_p[\zeta][z_1,\ldots,z_n]} and the reduction \m{\overline f} of \m{f} modulo \m{\maxideal} is non-zero. For all \m{x \in D} and all \m{\lambda \in \Lambda_x},
\algn{f(\zeta^{(b_{x})_1+\lambda x_1},\ldots,\zeta^{(b_{x})_n+\lambda x_n})&=f(\zeta^{(b_{x})_1}(\zeta^\lambda)^{ [x_1]},\ldots,\zeta^{(b_{x})_n}(\zeta^\lambda)^{ [x_n]})\\&=g_x(\zeta^{\lambda})=0,}
where, for \m{y \in R}, \m{[y]} is the integer in \m{\{0,\ldots,q-1\}} congruent to \m{y} modulo \m{q}, and \m{g_x \in \ZZ_p[\zeta][z]} is some polynomial depending on \m{x}. So \m{\prod_{\lambda \in \Lambda_x}(z-\zeta^{\lambda}) \mid g_x}, implying that  \m{(z-1)^{\lceil\varepsilon q\rceil} \mid \overline g_x}, in turn implying that
\algn{\TADICVALUATION\left(\overline g_x(1+t)\right)=\TADICVALUATION\left(\overline f((1+t)^{[x_1]},\ldots,(1+t)^{[x_n]})\right) \GEQSLANT \varepsilon q.} Let 
\algn{C_D = \left\{((1+t)^{[x_1]},\ldots,(1+t)^{[x_n]}) \subseteq \overline C^n\suchthat (x_1,\ldots,x_n)\in D\right\}.} Then, for all \m{s \in C_D}, \m{\TADICVALUATION\left(\overline f(s)\right)\GEQSLANT \varepsilon q}. As \m{\overline f \in \overline T[z_1,\ldots,z_n]} is non-zero, has total degree less than \m{\varepsilon \delta q/pnk\LEQSLANT q}, and has coefficients in \m{\FF}, it follows from lemma~\ref{vanlem} with some \m{c \in \FF^\times} and \m{\nu = \varepsilon/n} that \algn{\delta q^n \LEQSLANT |D| = |C_D| < pkq^{n-1} (\varepsilon \delta q/pnk) (n/\varepsilon) = \delta q^n,} which is a contradiction.
\end{proof}

\bibliographystyle{amsplain}

{

}

\end{document}